%% file: paper.tex
\documentclass{article}
\usepackage{graphicx}
\usepackage{amssymb}
\usepackage{amsmath}
\usepackage{amsthm}
\usepackage{color}

\newtheorem{lem}{Lemma}

\newtheorem{thm}{Theorem}

\newtheorem{pro}{Proposition}
\newtheorem{cor}{Corollary}

\def\bibitemart#1#2#3#4#5#6#7{\bibitem{#1}#2. #3. #4, #5:#7, #6.}

\begin{document}

\title{The equidistant dimension of some graphs of convex polytopes}
\author{Aleksandar Savi\'{c}\\
Faculty of Mathematics, University of Belgrade\\
Studentski trg 16/IV, 11 000 Belgrade, Serbia\\
\and
Zoran Maksimovi\'{c}\\
Military Academy, University of Defence\\
generala Pavla Juri\v{s}i\'{c}a \v{S}turma 33, 11 000 Belgrade, Serbia\\
\and
Milena Bogdanovi\'{c}\\
Faculty of Information Technologies, Belgrade \\
Metropolitan University,\\
Tadeu\v{s}a Ko\v{s}\'cu\v{s}ka 63, 11 000 Belgrade, Serbia\\
\and
Jozef Kratica\\
Mathematical Institute, Serbian Academy of Sciences and Arts\\
Kneza Mihaila 36/III, 11 000 Belgrade, Serbia
}

\maketitle

\begin{abstract}
This paper is devoted to some rotationally symmetric classes of graphs denoted in literature as convex polytope graphs. Exact value of equidistant dimension is found for $T_n$. Next, for even $n$ exact values are found for $R''_n$ and $S''_n$, while for odd $n$ exact value is found for $S_n$. Finally, for odd $n$, lower bound are found for $R''_n$ and $S''_n$. 
\end{abstract}


\section{Introduction}

Graph theory is mathematical discipline which concerns graphs and many aspects of their structre and applicability. One of more important aspects which helps in this study are graph invariants. From middle seventies was introduced one new graph invariant, notably \textit{metric dimension}. The metric dimension was independently introduced by Harary and Melter in \cite{har76} and by Slater in \cite{sla75}. The importance of metric dimesion and its variants  mirrors in great number of practical applications including location problems in networks of different nature (see \cite{cac07}). For example, in a computer network, which can be modeled as a graph, locating failure can be of great interest.  To do so, we can consider a subset of vertices $S$ in order that any vertex of the graph can be uniquely identified by its vector of distances to the vertices of $S$. This resolvability lies in the fundamental concept of metric dimension. Metric dimensions are concerned with the notion of these resolving sets, and especially such that their cardinality is minimal.

However, in recent times question of anonimity, especially in social networks, has lead to question where such unique identification is not possible as can be seen in papers \cite{cacsrom17,ches11,ches13,fed08,tru16,zhe11,zho08}. In answer to this tendencies was introduced one new metric dimension graph invariant, notably \text{equidistant metric dimension}, by Gonz\'{a}les et al. in \cite{gon16}. This new metric dimension can be used in problems of network's anonimization, as well as application of resolving sets in number theory. 

All graphs considered in this paper are connected, undirected, simple, and finite. Also, we define \textit{distance} between vertices $u$ and $v$, in notation $d(u,v)$ as number of edges in the shortest path between vertices $u$ and $v$. A vertex $x$ is notated as \textit{equidistant} from pair of vertices $u$ and $v$ if $d(u,x)=d(v,x)$. A subset $S$ of vertices is a \textit{distance-equalizer set} if for each pair of vertices $u , v \in V(G) \setminus S$ there exists a vertex $x \in S$ that is equidistant from $u$
and $v$. The \textit{equidistant dimension} of $G$, denoted $eqdim(G)$, is the minimum cardinality
of a distance-equalizer set of $G$, where \textit{equidistant base} of $G$ is distance equilizer set of cardinality $eqdim(G)$. Gispert-Fern{\'a}ndez and Rodr{\i}guez-Vel{\'a}zquez in \cite {np} proved that the equidistant dimension problem is NP-hard in a general case, and considered equidistant dimension of lexicographic product of graphs.

\subsection{Basic notions}

Standard notations for the set of vertices $V(G)$, the set of edges $E(G)$ will be used. The order of a graph $G$ is $|V(G)|$. Open neighborhood of a vertex $v$ from $G$ is the $N(v)=\{w\in V(G)|vw\in E(G)\}$, while its closed neighborhood is defined as $N[v]=N(v)\cup\{v\}$. By $deg(v)=|N(v)|$ we denote the degree of a vertex.

For a graph $G$, $\Delta(G)$ is its maximum degree $\Delta(G)=\max\{deg(v)|v\in V(G)\}$, $\delta(G)$ is its minimum degree $\delta(G)=\min\{deg(v)|v\in V(G)\}$. The diameter of a graph $G$ is the maximal possible distance between to vertices from the graph, i.e. $Diam(G)=\max\{d(u,v)|u,v\in V(G)\}$.

A clique is subset of pairwise adjacent vertices while independent set is subset of pairwise non-adjacent vertices. The clique number $\omega(G)$ and the independence number $\alpha(G)$ is the size of a largest clique/independent set of maximum cardinality, respectively.

For a pair of vertices $u,v\in V(G)$, as in \cite{bal09} we define the set 
$uWv$ of vertices that are equally distant from vertices $u$ and $v$, i.e.
$uWv=\{w\in V(G)|d(u,w)=d(v,w)\}$.

In this paper we will make important use of following statements:

\begin{thm} \label{ext1} \mbox{\rm(\cite{gon16})} For every graph $G$ of order $n \ge 2$, the following statements hold.
\begin{itemize}
\item $eqdim(G) = 1$ if and only if $\Delta(G) = |V| - 1$;
\item $eqdim(G) = 2$ if and only if $\Delta(G) = |V| - 2$.
\end{itemize}
\end{thm} 

\begin{cor} \label{ext1a} \mbox{\rm(\cite{gon16})}
If $G$ is a graph with $\Delta(G) < |V| - 2$ then $eqdim(G) \ge 3$.
\end{cor}

\begin{thm} \label{ext2} \mbox{\rm(\cite{gon16})} For every graph $G$, the following statements hold.
\begin{itemize}
\item If $n \ge 2$, then $eqdim(G) = |V| - 1$ if and only if $G$ is a path of order 2;
\item If $n \ge 3$, then $eqdim(G) = |V| - 2$ if and only if $G \in \{P_3, P_4, P_5, P_6,C_3,C_4,C_5\}$.
\end{itemize}
\end{thm} 

\begin{cor} \label{ext2a} \mbox{\rm(\cite{gon16})}
If $G$ is a graph of order $|V| \ge 7$, then $1 \le eqdim(G) \le |V|-3$.
\end{cor} 

\begin{pro} \label{ub1} \mbox{\rm(\cite{gon16})} For every graph $G$ of order $|V| \ge 2$, the following statements hold.
\begin{enumerate}
\item[I)] $eqdim(G) \le |V| - \Delta(G)$ and the bound is tight if $\Delta(G) \ge \frac{|V|}{2}$;
\item[II)] $eqdim(G) \le |V| - \omega(G)$, if $G \not\cong K_{|V|}$, and the bound is tight if $\omega(G) \le \frac{|V|}{2}$;
\item[III)] $eqdim(G) \le \frac{|V| \cdot (Diam(G)-1)+1}{Diam(G)}$ , and the bound is tight if $Diam(G) = 2$;
\item[IV)] $eqdim(G) \le |V| - \alpha(G)$, whenever $Diam(G) = 2$, and the bound is tight if $\alpha(G) \ge \frac{|V|}{2}$.
\end{enumerate}
\end{pro}

\begin{lem} \label{HSP} \mbox{\rm(\cite{VV})} Let $G$ be a graph, and $u$ and $v$ any vertices from $V(G)$. If $S$ is a distance-equalizer set of $G$, then 
	$S \, \bigcap \, (\{u,v\} \, \bigcup \, _u{W_v}) \neq \emptyset$.
\end{lem}

\begin{cor} \label{HSC} \mbox{\rm(\cite{VV})} Let $G$ be a graph, and $u$ and $v$ any vertices from $V(G)$.
	If $S$ is a distance-equalizer set of $G$ and $_u{W_v} = \emptyset$ then $u \in S$ or $v \in S$.
\end{cor}

\section{Main results}

A convex polytope is a polytope that is also a convex set of points in the $n$-dimensional Euclidean space $\mathbb{R}^n$. The graph of a convex polytope is constructed by the vertices and edges with the same incidence relation, for the ﬁrst time considered in \cite{bacya92}. These graphs are rotationally symmetric and planar graphs. Certain graph invariant of convex polytopes were
considered in \cite{kra12a}. All graphs of convex polytopes are defined by 2 cycles of $n$ vertices which build one inner and one outer face. Between these cycles are arrays of $n$ $m$-sided faces, with various $m$, determined by arrays of $n$ vertices, so order of this graphs is always multiple of $n$. For all these graphs it will be considered that $n \ge 5$, and it should be noted that vertex indices are taken modulo $n$.

\subsection{The equidistant dimension of $R''_n$}

\input{r2n.tex}

One class of convex polytopes that was introduced in \cite{mil06}, is $R''_n$. It consists of $2n$ 5-sided faces, $n$ 6-sided faces, and a pair of $n$-sided faces. They are defined with sets of vertices and edges. Formally,

$$V(R''_n)=\{ a_i, b_i, c_i, d_i, e_i, f_i \,|\, i=0,\ldots, n-1\},$$
and
$$E(R''_n)=\{ a_i a_{i+1}, f_i f_{i+1}, a_i b_i, c_i d_i, e_i f_i, b_i c_i, b_{i+1} c_{i},\, d_i e_i, d_{i+1} e_i|\, i=0,\ldots, n-1\}.$$

The tight lower bound of equidistant dimension of $R''_n$ is determined
in the following lemma.

\begin{lem} \label{l2} $eqdim(R''_n) \ge 3n$ \end{lem}
\begin{proof} 
Let $S$ be an equilizer set of $R''_n$.
Let us find $_{c_i}W_{d_i}=\{x\,|\,d(c_i,x)=d(d_i,x)\}$. It can be seen that for all $i$ and $j$ it holds 
$d(c_i,f_j) \neq d(d_i,f_j)$,
$d(c_i,e_j) \neq d(d_i,e_j)$,
$d(c_i,a_j) \neq  d(d_i,a_j)$,
$d(c_i,b_j) \neq d(d_i,b_j)$,
$d(c_i,c_j) \neq d(d_i,c_j)$ and
$d(c_i,d_j) \neq d(d_i,d_j)$.
Therefore $_{c_i}W_{d_i}=\emptyset$. Corolary \ref{HSC} implies that either $c_i\in S$ or $d_i\in S$.

Next, let us find $_{b_i}W_{e_{i-1}}=\{x\,|\,d(b_i,x)=d(e_{i-1},x)\}$. Also for all $i$ and $j$ it holds
$d(b_i,a_j) \neq d(e_{i-1},a_j)$,
$d(b_i,b_j) \neq d(e_{i-1},b_j)$,
$d(b_i,e_j) \neq d(e_{i-1},e_j)$, 
$d(b_i,f_j) \neq d(e_{i-1},f_j)$,
$d(b_i,c_j) \neq d(e_{i-1},c_j)$ and 
$d(b_i,d_j) \neq d(e_{i-1},d_j)$.
Therefore $_{b_i}W_{e_{i-1}}=\emptyset$. Corolary \ref{HSC} implies that either $b_i\in S$ or $e_{i-1} \in S$.

Finally, let us find $_{a_i}W_{f_{i-1}}=\{x\,|\,d(a_i,x)=d(f_{i-1},x)\}$. Also for all $i$ and $j$ it holds
$d(a_i,a_j) \neq d(f_{i-1},a_j)$,
$d(a_i,b_j) \neq d(f_{i-1},b_j)$,
$d(a_i,e_j) \neq d(f_{i-1},e_j)$, 
$d(a_i,f_j) \neq d(f_{i-1},f_j)$,
$d(a_i,c_j) \neq d(f_{i-1},c_j)$ and 
$d(a_i,d_j) \neq d(f_{i-1},d_j)$.
Therefore $_{a_i}W_{f_{i-1}}=\emptyset$. Corolary \ref{HSC} implies that either $a_i\in S$ or $f_{i-1} \in S$.

Since we have $3n$ disjoint pairs of vertices, with at least one vertex from each pair including in $S$, 
so $|S|\ge 3n$.
\end{proof}

The equidistant dimension of $R''_n$ for even $n$ is determined
in the following theorem.

\begin{thm} For even $n$, $eqdim(R''_n)=3n$. \end{thm}
\begin{proof}

Lower bound directly follows from Lemma \ref{l2}. 
Let $S=\{a_i,b_i,c_i\, |\,0\le i\le n-1 \}$. The set $S$ is a distance equalizer set for $R''_n$ which follows from the fact that for every pair of vertices $u, v$ from $V(R''_n) \setminus S$ there is a vertex $x$ from $S$ that is equidistant from $u$ and $v$, i.e. $d(u,x) = d(v,x)$. The equidistant vertices
are noted in table \ref{tabr2even}.
\end{proof}

\begin{table}
\caption{Distance-equalizer of $R''_n$, when $n$ is even} \label{tabr2even}
 \small
\begin{center} 
\begin{tabular}{|c|c|c|c|}
\hline
$u$ & $v$ &  cond. & $x$ \\
\hline
$d_0$ & $d_i$ & $i$ odd &  $b_\frac{i+1}{2}$  \\
  &  & $i$ even &   $c_\frac{i}{2}$\\
 \hline
 $d_0$ & $e_0$ &  & $a_2$   \\
  $d_0$ & $e_{n-1}$ &  & $c_3$   \\
$d_0$ & $e_i$ & other odd $i$ &   $a_{\frac{n+i+1}{2}}$ \\
&  & other even $i$  &  $a_{\frac{i}{2}+1}$ \\
\hline
$d_0$ & $f_i$ & $i=0$ or $i=1$ &  $b_2$ \\
$d_0$ & $f_{n-1}$ &  &  $b_3$ \\
$d_0$ & $f_{n-2}$ &  &  $c_3$ \\

$d_0$ & $f_i$ & other odd  $i$  &  $a_{\frac{i+3}{2}} $ \\

&  & other even $i$ &   $a_{\frac{n+i}{2}} $ \\
\hline
$e_0$ & $e_i$ & $i$ odd &  $c_\frac{i+1}{2}$ \\
&  & $i$ even &  $b_{\frac{i}{2}+1}$ \\
\hline
$e_0$ & $f_0$ &  &  $c_2$ \\
$e_0$ & $f_1$ & &  $a_3$ \\
$e_0$ & $f_{n-4}$ &  &  $b_{n-1}$ \\
$e_0$ & $f_i$ & $n-3\le i \le n-1$&  $a_{n-1}$ \\
$e_0$ & $f_i$ & other odd $i$&  $a_{\frac{i+3}{2}}$ \\
$e_0$ & $f_i$ & other even $i$&  $c_{\frac{n+i}{2}+1}$ \\
\hline
$f_0$ & $f_i$ & $i$ odd &  $c_\frac{i+1}{2}$ \\
	&  & $i$ even &  $b_{\frac{i}{2}+1}$ \\
	\hline
\end{tabular}
\end{center}
\end{table}

\subsection{The equidistant dimension of $S_n$}

This class of convex polytopes $S_n$ was introduced in \cite{imr11c} where authors
also determined its metric dimension. It consists of $2n$ 3-sided faces, $2n$ 4-sided faces, and a pair of $n$-sided faces. Mathematically, the set of vertices and the set of edges are 

$$V(S_n)=\{ a_i, b_i, c_i, d_i \,|\, i=0,\ldots, n-1\},$$
$$E(S_n)=\{ a_i a_{i+1}, b_i b_{i+1}, c_i c_{i+1}, d_i d_{i+1},
a_i b_i, b_i c_i, c_i d_i, a_{i+1} b_i \,|\, i=0,\ldots, n-1\}$$

\input{sn.tex}

The following theorem determines equidistant dimension of $S_n$, for odd $n$.

\begin{thm} For odd $n$, it holds $eqdim(S_n)=2n$.\end{thm}
\begin{proof}

Step 1: For odd $n$, $eqdim(S_n) \ge 2n$ \\
Let $S$ be an equilizer set of $S_n$.
Let us find $_{a_i}W_{b_{i+2}}=\{x\,|\,d(a_i,x)=d(b_{i+2},x)\}$. It can be seen that for all $i$ and $j$ it holds 
$d(a_i,a_j) \neq d(b_{i+2},a_j)$,
$d(a_i,b_j) \neq d(b_{i+2},b_j)$,
$d(a_i,c_j) \neq d(b_{i+2},c_j)$,
$d(a_i,d_j) \neq d(b_{i+2},d_j)$,
Therefore $_{a_i}W_{b_{i+2}}=\emptyset$. Corolary \ref{HSC} implies that either $a_i\in S$ or $b_{i+2}\in S$.

Next, let us find $_{c_i}W_{d_{i}}=\{x\,|\,d(c_i,x)=d(d_{i},x)\}$. It is easy to see that if $x=d_j$ for each $j$ it holds
$d(c_i,x)=d(d_{i},x)+1$. Similarly, if $x \in \{a_j,b_j,c_j | 0 \le j \le n-1\}$ it holds $d(c_i,x)=d(d_{i},x)-1$. 
Since  $_{c_i}W_{d_{i}}=\emptyset$, then by Corolary \ref{HSC} it holds that either $c_i\in S$ or $d_{i} \in S$.

Since we have $2n$ disjoint pairs of vertices, with at least one vertex from each pair including in $S$, 
so $|S|\ge 2n$.

Step 2: For odd $n$, $eqdim(S_n) \le 2n$ \\

Let $S= \{a_i,c_i | 0 \le i \le n-1\}$. In Table \ref{tabsn} are presented all pairs of vertices from $V(S_n) \setminus S$, $n$ odd, and are shown their respective equidistant ones.
As can be seen from Table \ref{tabsn} for each pair $u$ and $v$ of vertices there exist at least one vertex $x \in S$ equidistant from them, 
i.e. $d(u,x)=d(v,x)$. Therefore, the set $S$ is a distance equalizer set for $S_n$, so $eqdim(S_n)=2n$, for odd $n$.

\end{proof}

\begin{table}
\caption{Distance-equalizer of $S_n$, when $n$ is odd} \label{tabsn}
 \small
\begin{center} 
\begin{tabular}{|c|c|c|c|}
\hline
$u$ & $v$ &  cond. & $x$ \\
\hline
$b_0$ & $b_i$ & $i$ odd &  $a_\frac{i+1}{2}$  \\
  &  & $i$ even &   $c_\frac{i}{2}$\\
 \hline
 $b_0$ & $d_i$ & $i$ odd & $c_\frac{n+i}{2}$   \\
  &  & $i$ even & $c_\frac{i}{2}$   \\
\hline
$d_0$ & $d_i$ & $i$ odd &  $c_\frac{n+i}{2}$ \\
	&  & $i$ even &  $c_\frac{i}{2}$ \\
	\hline
\end{tabular}
\end{center}
\end{table}

\subsection{The equidistant dimension of $S''_n$}

Another class of convex polytopes $S''_n$ was introduced in \cite{imr11b}. Authors
also determined its metric dimension in mentioned paper. Similar to $S_n$ it also consists of $2n$ 3-sided faces, $2n$ 4-sided faces, and a pair of $n$-sided faces. Set of vertices and set od edges that are
defining these convex polytopes are

$$V(S''_n)=\{ a_i, b_i, c_i, d_i \},$$ 
$$E(S''_n)=\{ a_i a_{i+1}, b_i b_{i+1}, c_i c_{i+1}, d_i d_{i+1}, a_i b_i, b_i c_i, c_i d_i, b_{i+1} c_i \}$$ 

 
\input{s2n.tex}

The tight lower bound of equidistant dimension of $S''_n$ is determined
in the following lemma.

\begin{lem}\label{l1} $eqdim(S''_n) \ge 2n$.\end{lem}
\begin{proof}
Let $S$ is a equilizer set of $S''_n$.
Let us find $_{a_i}W_{b_{i}}=\{x\,|\,d(a_i,x)=d(b_{i},x)\}$. It can be seen that for all $i$ and $j$ it holds 
$d(a_i,a_j) \neq d(b_{i},a_j)$,
$d(a_i,b_j) \neq d(b_{i},b_j)$,
$d(a_i,c_j) \neq d(b_{i},c_j)$,
$d(a_i,d_j) \neq d(b_{i},d_j)$,
Therefore $_{a_i}W_{b_{i}}=\emptyset$. Corolary \ref{HSC} implies that either $a_i\in S$ or $b_{i}\in S$.

Next, let us find $_{c_i}W_{d_{i}}=\{x\,|\,d(c_i,x)=d(d_{i},x)\}$. It is easy to see 
$d(c_i,a_j) \neq d(d_{i},a_j)$,
$d(c_i,b_j) \neq d(d_{i},b_j)$,
$d(c_i,c_j) \neq d(d_{i},c_j)$,
$d(c_i,d_j) \neq d(d_{i},d_j)$,
Since  $_{c_i}W_{d_{i}}=\emptyset$, then by Corolary \ref{HSC} it holds that either $c_i\in S$ or $d_{i} \in S$.

Since we have $2n$ disjoint pairs of vertices, with at least one vertex from each pair including in $S$, 
so $|S|\ge 2n$.

\end{proof}

The following theorem determines eqidistant dimension of $S''_n$, for even $n$.

\begin{thm} For even $n$, it holds $eqdim(S''_n)=2n$.\end{thm}
\begin{proof}

Lower bound directly follows from Lemma \ref{l1}. Let $S= \{a_i,c_i | 0 \le i \le n-1\}$. 

In Table \ref{tabs2n} are presented all pairs of vertices from $V(S''_n) \setminus S$, $n$ even, and are shown their respective equidistant ones.

As can be seen from Table \ref{tabs2n} for each pair $u$ and $v$ there exist at least one vertex $x \in S$ equidistant from them, 
i.e. $d(u,x)=d(v,x)$. Therefore, the set $S$ is a distance equalizer set for $S''_n$, so $eqdim(S''_n)=2n$, for even $n$.

\end{proof}

\begin{table}
\caption{Distance-equalizer of $S''_n$, when $n$ is even} \label{tabs2n}
 \small
\begin{center} 
\begin{tabular}{|c|c|c|c|}
\hline
$u$ & $v$ &  cond. & $x$ \\
\hline
$b_0$ & $b_i$ & $i$ odd &  $c_\frac{i-1}{2}$  \\
  &  & $i$ even &   $a_\frac{i}{2}$\\
 \hline
 $b_0$ & $d_i$ & $i=n-1$ & $c_{n-1}$ \\
 &  & other odd $i$ & $a_\frac{n+i-1}{2}$   \\
  &  & $i$ even & $c_\frac{i}{2}$   \\
\hline
$d_0$ & $d_i$ & $i$ odd &  $a_\frac{i+1}{2}$ \\
	&  & $i$ even &  $c_\frac{i}{2}$ \\
	\hline
\end{tabular}
\end{center}
\end{table}

\subsection{The equidistant dimension of $T_n$}

Yet another class of convex polytopes $T_n$ was introduced in \cite{imr11c}. In the same paper, authors
determined its metric dimension. It consists of $4n$ 3-sided faces, $n$ 4-sided faces, and a pair of $n$-sided faces. Formally, set of vertices and set od edges are

$$V(T_n)=\{ a_i, b_i, c_i, d_i \,|\, i=0,\ldots, n-1\},$$
$$E(T_n)=\{ a_i a_{i+1}, b_i b_{i+1}, c_i c_{i+1}, d_i d_{i+1},
a_i b_i, b_i c_i, c_i d_i, a_{i+1} b_i, c_{i} d_{i+1} \,|\, i=0,\ldots, n-1\}$$

\input{tn.tex}

\begin{table}
\caption{Distance-equalizer of $T_n$} \label{eqt}
 \small
\begin{center} 
\begin{tabular}{|c|c|c|c|}
\hline
$u$ & $v$ &  cond. & $x$ \\
\hline
$c_0$ & $c_i$ & $i$ odd &  $a_\frac{i+1}{2}$ \\
&  & $i$ even &  $b_\frac{i}{2}$ \\
\hline
$c_0$ & $d_i$ & $i$ odd & $b_\frac{i+1}{2}$  \\
&  & $i$ even, $i\ge2$ & $a_{\frac{i}{2}+1}$  \\
&  & $i=0$ & $b_1$\\
\hline
$d_0$ & $d_i$ & $i$ odd &  $b_\frac{i+1}{2}$ \\
  &  & $i$ even &  $a_{\frac{i}{2}+1}$ \\
 \hline
\end{tabular}
\end{center}
\end{table}    

The following theorem determines $eqdim(T_n)$.

\begin{thm}$eqdim(T_n)=2n$.\end{thm}
\begin{proof}
Step 1: $eqdim(T_n) \ge 2n$ \\
Let $S$ is a equilizer set of $T_n$.
Let us find $_{b_i}W_{c_i}=\{x\,|\,d(b_i,x)=d(c_i,x)\}$. It is obvious that for all $i$ and $j$ it holds 
$d(b_i,a_j)=d(c_i,a_j)-1$,
$d(b_i,b_j)=d(c_i,b_j)-1$,
$d(b_i,c_j)=d(c_i,c_j)+1$ and
$d(b_i,d_j)=d(c_i,d_j)+1$.
Therefore $_{b_i}W_{c_i}=\emptyset$. Corolary \ref{HSC} implies that either $b_i\in S$ or $c_i\in S$.

Finally, let us find $_{a_i}W_{d_{i-1}}=\{x\,|\,d(a_i,x)=d(d_{i-1},x)\}$. It is obvious that for all $i$ and $j$ it holds
$d(a_i,a_j)=d(d_{i-1},a_j)-X$,
$d(a_i,b_j)=d(d_{i-1},b_j)-X$,
$d(a_i,c_j)=d(d_{i-1},c_j)+X$ and 
$d(a_i,d_j)=d(d_{i-1},d_j)+X$.
Therefore $_{a_i}W_{d_{i-1}}=\emptyset$. Corolary \ref{HSC} implies that either $a_i\in S$ or $d_{i-1} \in S$.

Since we have $2n$ disjoint pairs of vertices, with at least one vertex from each pair including in $S$, 
so $|S|\ge 2n$.

Step 2: $eqdim(T_n) \le 2n$ \\
Let $S=\{a_i,b_i\, |\,0\le i\le n-1 \}$. 

In Table \ref{eqt} are presented all pairs of vertices from $V(T_n) \setminus S$, and are shown their respective equidistant ones.

As can be seen from Table \ref{eqt} for each pair $u$ and $v$ there exist at least one vertex $x \in S$ equidistant from them, 
i.e. $d(u,x)=d(v,x)$. Therefore, the set $S$ is a distance equalizer set for $T_n$, so $eqdim(T_n)=2n$.

\end{proof}

\section{Discussion}

Metric dimension of graphs of convex polytopes $R''_n$, $S_n$, $S''_n$ and $T_n$ is constant and equal to 3, see \cite{imr11b,imr11c,imr11e}. 

However, their equidistant dimension depends on $n$: 
\begin{itemize}
\item For even $n$, $eqdim(R''_n) = 3n$ and $eqdim(S''_n) = 2n$;
\item For odd $n$, $eqdim(S_n) = 2n$;
\item For each $n$, $eqdim(T_n) = 2n$.
\end{itemize}

Looking at the literature it can be seen that doubly and strong metric dimension have also quite different exact values than equidistant dimension. For example, in \cite{kra12a} is shown that $T_n$ has doubly metric dimension equal to 3 for $n \neq 7$ and strong metric dimension of $\frac{5n}{2}$ for even $n$. Nevertheless, for odd $n$ strong metric dimension is equal to equidistant dimension ($2n$), but it is purely coincidental, since strong metric base is $\{a_i,d_i | 0 \le i \le n-1\}$ while equidistant base is $\{a_i,b_i | 0 \le i \le n-1\}$.

It should be interesting to compare these values with lower and upper bounds from literature.
It is obvious that lower bound from Corollary \ref{ext1a} is equal to 3 which is much smaller than values found in this paper. Upper bound from Corollary \ref{ext2a} is equal to $|V(G)| -2$ is much larger than mentioned values. 

Next, it holds $\Delta(R''_n)=3$ and $\Delta(S_n)=\Delta(S''_n)=\Delta(T_n)=5$, so by Propostion \ref{ub1} part I) is also much larger than mentioned values. Since $\omega(R''_n)=2$ and $\omega(S_n)=\omega(S''_n)=\omega(T_n)=3$ then Propostion \ref{ub1} Part II) 
produces similar results. Computation in Propostion \ref{ub1} Part III) is more complex, but final bounds are also close to $|V(G)|$. 
The upper bound from Propostion \ref{ub1} Part IV) is not close to $|V(G)|$ but also far from mentioned values.

\section{Conclusions}

In this paper equidistant dimension of some convex polytopes are considered.
Exact value is equal to $2n$ for $T_n$, $S_n$ with odd $n$ and $S''_n$ with even $n$. For 
$R''_n$ with even $n$ exact value is equal to $3n$. Finally, for odd $n$, lower bounds of $3n$ and $2n$ are found for $R''_n$ and $S''_n$, respectively.

Future work can be directed towards obtaining equidistant dimension of
some other classes of graphs, and similar graph invariants for convex polytopes.
Also, it would be of interest to solve equidistant dimension problem by exact method
or heuristic.

\end{document}

%% file: r2n.tex
\begin{figure}[htbp]
    \centering\setlength\unitlength{1mm}
\setlength\unitlength{1mm}
\begin{picture}(130,65)
\thicklines
\tiny
\put(47.7,5.0){\circle*{1}} \put(47.7,5.0){\line(1,1){7.3}} \put(47.7,5.0){\line(-5,3){8.7}} 
\put(55.0,12.3){\circle*{1}} \put(55.0,12.3){\line(4,1){10.0}} \put(55.0,12.3){\line(-3,5){5.0}} 
\put(65.0,15.0){\circle*{1}} \put(65.0,15.0){\line(4,-1){10.0}} \put(65.0,15.0){\line(0,1){10.0}} 
\put(75.0,12.3){\circle*{1}} \put(75.0,12.3){\line(1,-1){7.3}} \put(75.0,12.3){\line(3,5){5.0}} 
\put(82.3,5.0){\circle*{1}} \put(82.3,5.0){\line(5,3){8.7}} 
\put(39.0,10.0){\circle*{1}} \put(39.0,10.0){\line(-5,-2){12.7}} \put(39.0,10.0){\line(-1,6){2.3}} 
\put(50.0,21.0){\circle*{1}} \put(50.0,21.0){\line(-6,1){13.3}} \put(50.0,21.0){\line(2,5){4.6}} 
\put(65.0,25.0){\circle*{1}} \put(65.0,25.0){\line(-5,4){10.4}} \put(65.0,25.0){\line(5,4){10.4}} 
\put(80.0,21.0){\circle*{1}} \put(80.0,21.0){\line(-2,5){4.6}} \put(80.0,21.0){\line(6,1){13.3}} 
\put(91.0,10.0){\circle*{1}} \put(91.0,10.0){\line(1,6){2.3}} \put(91.0,10.0){\line(5,-2){12.7}} 
\put(26.4,5.4){\circle*{1}} \put(26.4,5.4){\line(-4,1){9.7}} 
\put(36.7,23.3){\circle*{1}} \put(36.7,23.3){\line(-1,1){7.1}} 
\put(54.6,33.6){\circle*{1}} \put(54.6,33.6){\line(-1,4){2.6}} 
\put(75.4,33.6){\circle*{1}} \put(75.4,33.6){\line(1,4){2.6}} 
\put(93.3,23.3){\circle*{1}} \put(93.3,23.3){\line(1,1){7.1}} 
\put(103.6,5.4){\circle*{1}} \put(103.6,5.4){\line(4,1){9.7}} 
\put(16.7,7.9){\circle*{1}} \put(16.7,7.9){\line(0,1){14.6}} 
\put(29.6,30.4){\circle*{1}} \put(29.6,30.4){\line(-3,-2){12.3}} \put(29.6,30.4){\line(2,3){7.9}} 
\put(52.1,43.3){\circle*{1}} \put(52.1,43.3){\line(-1,0){14.6}} \put(52.1,43.3){\line(2,1){12.9}} 
\put(77.9,43.3){\circle*{1}} \put(77.9,43.3){\line(-2,1){12.9}} \put(77.9,43.3){\line(1,0){14.6}} 
\put(100.4,30.4){\circle*{1}} \put(100.4,30.4){\line(-2,3){7.9}} \put(100.4,30.4){\line(3,-2){12.3}} 
\put(113.3,7.9){\circle*{1}} \put(113.3,7.9){\line(0,1){14.6}} 
\put(17.4,22.5){\circle*{1}} \put(17.4,22.5){\line(-5,3){8.7}} 
\put(37.5,42.6){\circle*{1}} \put(37.5,42.6){\line(-3,5){5.0}} 
\put(65.0,50.0){\circle*{1}} \put(65.0,50.0){\line(0,1){10.0}} 
\put(92.5,42.6){\circle*{1}} \put(92.5,42.6){\line(3,5){5.0}} 
\put(112.6,22.5){\circle*{1}} \put(112.6,22.5){\line(5,3){8.7}} 
\put(8.7,27.5){\circle*{1}} \put(8.7,27.5){\line(1,1){23.8}} 
\put(32.5,51.3){\circle*{1}} \put(32.5,51.3){\line(4,1){32.5}} 
\put(65.0,60.0){\circle*{1}} \put(65.0,60.0){\line(4,-1){32.5}} 
\put(97.5,51.3){\circle*{1}} \put(97.5,51.3){\line(1,-1){23.8}} 
\put(121.3,27.5){\circle*{1}} 
\thinlines
\put(47.7,5.0){\line(-1,-2){3}} \put(82.3,5.0){\line(1,-2){3}}
\put(26.4,5.4){\line(2,-3){3}} \put(103.6,5.4){\line(-2,-3){3}}
\put(16.7,7.9){\line(-2,-3){3}} \put(113.3,7.9){\line(2,-3){3}}
\put(8.7,27.5){\line(-1,-4){3}} \put(121.3,27.5){\line(1,-4){3}}
\put(49.1,3.0){$a_2$} \put(55.0,8.9){$a_1$} \put(63.0,11.0){$a_0$} \put(71.0,8.9){$a_{n-1}$} \put(76.9,3.0){$a_{n-2}$} 
\put(34.4,11.5){$b_2$} \put(46.5,23.6){$b_1$} \put(63.0,28.0){$b_0$} \put(79.5,23.6){$b_{n-1}$} \put(91.6,11.5){$b_{n-2}$} 
\put(23.9,7.5){$c_2$} \put(35.5,25.8){$c_1$} \put(54.7,36.0){$c_0$} \put(76.5,35.1){$c_{n-1}$} \put(94.8,23.5){$c_{n-2}$} \put(105.0,4.3){$c_{n-3}$} 
\put(12.8,8.5){$d_2$} \put(26.2,31.8){$d_1$} \put(49.5,45.2){$d_0$} \put(76.5,45.2){$d_{n-1}$} \put(99.8,31.8){$d_{n-2}$} \put(114,8.5){$d_{n-3}$} 
\put(15.8,25.4){$e_1$} \put(37.5,45.4){$e_0$} \put(66.2,52.0){$e_{n-1}$} \put(94.4,43.2){$e_{n-2}$} \put(114.4,21.5){$e_{n-3}$}
\put(6.0,28.5){$f_1$} \put(30.5,53.0){$f_0$} \put(64.0,62.0){$f_{n-1}$} \put(97.5,53.0){$f_{n-2}$} \put(122.0,28.5){$f_{n-3}$}

\end{picture}
    \caption{\label{fig:r2n} \small The graph of convex polytope $R''_n$}
\end{figure}

%% file: sn.tex
\begin{figure}[htbp]
    \centering\setlength\unitlength{1mm}

\setlength\unitlength{1mm}
\begin{picture}(130,45)
\thicklines
\tiny
\put(47.7,5.0){\circle*{1}} \put(47.7,5.0){\line(1,1){7.3}} \put(47.7,5.0){\line(-3,-1){8.8}} \put(47.7,5.0){\line(-1,5){1.8}} 
\put(55.0,12.3){\circle*{1}} \put(55.0,12.3){\line(4,1){10.0}} \put(55.0,12.3){\line(-5,1){9.1}} \put(55.0,12.3){\line(1,3){3.0}} 
\put(65.0,15.0){\circle*{1}} \put(65.0,15.0){\line(4,-1){10.0}} \put(65.0,15.0){\line(-5,4){7.0}} \put(65.0,15.0){\line(5,4){7.0}} 
\put(75.0,12.3){\circle*{1}} \put(75.0,12.3){\line(1,-1){7.3}} \put(75.0,12.3){\line(-1,3){3.0}} \put(75.0,12.3){\line(5,1){9.1}} 
\put(82.3,5.0){\circle*{1}} \put(82.3,5.0){\line(1,5){1.8}} \put(82.3,5.0){\line(3,-1){8.8}} 
\put(38.9,2.0){\circle*{1}} \put(38.9,2.0){\line(3,5){7.0}} \put(38.9,2.0){\line(-4,1){9.7}} 
\put(45.9,14.1){\circle*{1}} \put(45.9,14.1){\line(5,3){12.1}} \put(45.9,14.1){\line(-1,1){7.1}} 
\put(58.0,21.1){\circle*{1}} \put(58.0,21.1){\line(1,0){14.0}} \put(58.0,21.1){\line(-1,4){2.6}} 
\put(72.0,21.1){\circle*{1}} \put(72.0,21.1){\line(5,-3){12.1}} \put(72.0,21.1){\line(1,4){2.6}} 
\put(84.1,14.1){\circle*{1}} \put(84.1,14.1){\line(3,-5){7.0}} \put(84.1,14.1){\line(1,1){7.1}} 
\put(91.1,2.0){\circle*{1}} \put(91.1,2.0){\line(4,1){9.7}} 
\put(29.3,4.6){\circle*{1}} \put(29.3,4.6){\line(3,5){9.6}} \put(29.3,4.6){\line(-4,1){9.7}} 
\put(38.8,21.2){\circle*{1}} \put(38.8,21.2){\line(5,3){16.6}} \put(38.8,21.2){\line(-1,1){7.1}} 
\put(55.4,30.7){\circle*{1}} \put(55.4,30.7){\line(1,0){19.2}} \put(55.4,30.7){\line(-1,4){2.6}} 
\put(74.6,30.7){\circle*{1}} \put(74.6,30.7){\line(5,-3){16.6}} \put(74.6,30.7){\line(1,4){2.6}} 
\put(91.2,21.2){\circle*{1}} \put(91.2,21.2){\line(3,-5){9.6}} \put(91.2,21.2){\line(1,1){7.1}} 
\put(100.7,4.6){\circle*{1}} \put(100.7,4.6){\line(4,1){9.7}} 
\put(19.6,7.2){\circle*{1}} \put(19.6,7.2){\line(3,5){12.2}} 
\put(31.8,28.2){\circle*{1}} \put(31.8,28.2){\line(5,3){21.1}} 
\put(52.8,40.4){\circle*{1}} \put(52.8,40.4){\line(1,0){24.3}} 
\put(77.2,40.4){\circle*{1}} \put(77.2,40.4){\line(5,-3){21.1}} 
\put(98.2,28.2){\circle*{1}} \put(98.2,28.2){\line(3,-5){12.2}} 
\put(110.4,7.2){\circle*{1}} 
\thinlines
\put(47.7,5.0){\line(-1,-2){2}} \put(82.3,5.0){\line(1,-2){2}}
\put(38.9,2.0){\line(1,-1){4}} \put(91.1,2.0){\line(-1,-1){4}}
\put(38.9,2.0){\line(0,-1){5}} \put(91.1,2.0){\line(0,-1){5}}
\put(29.3,4.6){\line(0,-1){5}} \put(100.7,4.6){\line(0,-1){5}}
\put(19.6,7.2){\line(0,-1){5}} \put(110.4,7.2){\line(0,-1){5}}

\put(48.3,3.5){$a_2$} \put(54.5,9.7){$a_1$} \put(63.0,12.0){$a_0$} \put(71,9.7){$a_{n-1}$} \put(76.6,3.5){$a_{n-2}$} 
\put(37.6,4.3){$b_2$} \put(45.8,16.3){$b_1$} \put(58.9,22.5){$b_0$} \put(73.3,21.4){$b_{n-1}$} \put(85.3,13.2){$b_{n-2}$} \put(91.5,0.1){$b_{n-3}$} 
\put(27.9,6.7){$c_2$} \put(38.6,23.2){$c_1$} \put(56.1,32.2){$c_0$} \put(75.7,31.1){$c_{n-1}$} \put(92.2,20.4){$c_{n-2}$} \put(101.2,2.9){$c_{n-3}$} 
\put(16.4,7.7){$d_2$} \put(29.4,29.6){$d_1$} \put(51.3,42.3){$d_0$} \put(76.7,42.3){$d_{n-1}$} \put(98.6,29.6){$d_{n-2}$} \put(111.3,7.7){$d_{n-3}$} 

\end{picture}
    \caption{\label{fig:sn} \small The graph of convex polytope $S_n$}
\end{figure}

%% file: s2n.tex
\begin{figure}[htbp]
    \centering\setlength\unitlength{1mm}

\setlength\unitlength{1mm}
\begin{picture}(130,45)
\thicklines
\tiny
\put(47.7,5.0){\circle*{1}} \put(47.7,5.0){\line(1,1){7.3}} \put(47.7,5.0){\line(-5,3){8.7}} 
\put(55.0,12.3){\circle*{1}} \put(55.0,12.3){\line(4,1){10.0}} \put(55.0,12.3){\line(-3,5){5.0}} 
\put(65.0,15.0){\circle*{1}} \put(65.0,15.0){\line(4,-1){10.0}} \put(65.0,15.0){\line(0,1){10.0}} 
\put(75.0,12.3){\circle*{1}} \put(75.0,12.3){\line(1,-1){7.3}} \put(75.0,12.3){\line(3,5){5.0}} 
\put(82.3,5.0){\circle*{1}} \put(82.3,5.0){\line(5,3){8.7}} 
\put(39.0,10.0){\circle*{1}} \put(39.0,10.0){\line(1,1){11.0}} \put(39.0,10.0){\line(-5,-3){9.8}} \put(39.0,10.0){\line(0,1){11.2}} 
\put(50.0,21.0){\circle*{1}} \put(50.0,21.0){\line(4,1){15.0}} \put(50.0,21.0){\line(-1,0){11.2}} \put(50.0,21.0){\line(3,5){5.4}} 
\put(65.0,25.0){\circle*{1}} \put(65.0,25.0){\line(4,-1){15.0}} \put(65.0,25.0){\line(-5,3){9.6}} \put(65.0,25.0){\line(5,3){9.6}} 
\put(80.0,21.0){\circle*{1}} \put(80.0,21.0){\line(1,-1){11.0}} \put(80.0,21.0){\line(-3,5){5.4}} \put(80.0,21.0){\line(1,0){11.2}} 
\put(91.0,10.0){\circle*{1}} \put(91.0,10.0){\line(0,1){11.2}} \put(91.0,10.0){\line(5,-3){9.8}} 
\put(29.3,4.6){\circle*{1}} \put(29.3,4.6){\line(-4,1){9.7}} \put(29.3,4.6){\line(3,5){9.6}} 
\put(38.8,21.2){\circle*{1}} \put(38.8,21.2){\line(-1,1){7.1}} \put(38.8,21.2){\line(5,3){16.6}} 
\put(55.4,30.7){\circle*{1}} \put(55.4,30.7){\line(-1,4){2.6}} \put(55.4,30.7){\line(1,0){19.2}} 
\put(74.6,30.7){\circle*{1}} \put(74.6,30.7){\line(1,4){2.6}} \put(74.6,30.7){\line(5,-3){16.6}} 
\put(91.2,21.2){\circle*{1}} \put(91.2,21.2){\line(1,1){7.1}} \put(91.2,21.2){\line(3,-5){9.6}} 
\put(100.7,4.6){\circle*{1}} \put(100.7,4.6){\line(4,1){9.7}} 
\put(19.6,7.2){\circle*{1}} \put(19.6,7.2){\line(3,5){12.2}} 
\put(31.8,28.2){\circle*{1}} \put(31.8,28.2){\line(5,3){21.1}} 
\put(52.8,40.4){\circle*{1}} \put(52.8,40.4){\line(1,0){24.3}} 
\put(77.2,40.4){\circle*{1}} \put(77.2,40.4){\line(5,-3){21.1}} 
\put(98.2,28.2){\circle*{1}} \put(98.2,28.2){\line(3,-5){12.2}} 
\put(110.4,7.2){\circle*{1}} 
\thinlines
\put(47.7,5.0){\line(-1,-2){3}} \put(82.3,5.0){\line(1,-2){3}}
\put(39.0,10.0){\line(-1,-2){3}} \put(91.0,10.0){\line(1,-2){3}}
\put(29.3,4.6){\line(0,-1){6}} \put(100.7,4.6){\line(0,-1){6}}
\put(19.6,7.2){\line(0,-1){6}} \put(110.4,7.2){\line(0,-1){6}}

\put(48.3,3.5){$a_2$} \put(54.5,9.7){$a_1$} \put(63.0,12.0){$a_0$} \put(71,9.7){$a_{n-1}$} \put(76.5,3.5){$a_{n-2}$} 
\put(36.3,11.0){$b_2$} \put(48.0,22.7){$b_1$} \put(64.0,27.0){$b_0$} \put(80.0,22.7){$b_{n-1}$} \put(91.7,10.2){$b_{n-2}$} 
\put(27.9,6.7){$c_2$} \put(38.6,23.2){$c_1$} \put(56.1,32.2){$c_0$} \put(75.7,31.1){$c_{n-1}$} \put(92.2,20.4){$c_{n-2}$} \put(101.2,2.9){$c_{n-3}$} 
\put(16.3,7.7){$d_2$} \put(29,29.6){$d_1$} \put(51.3,42.3){$d_0$} \put(76.7,42.3){$d_{n-1}$} \put(98.6,29.6){$d_{n-2}$} \put(111.3,7.7){$d_{n-3}$} 

\end{picture}
    \caption{\label{fig:s2n} \small The graph of convex polytope $S''_n$}
\end{figure}

%% file: tn.tex
\begin{figure}[htbp]
    \centering\setlength\unitlength{1mm}

\setlength\unitlength{1mm}
\begin{picture}(130,50)
\thicklines
\tiny
\put(47.7,5.0){\circle*{1}} \put(47.7,5.0){\line(1,1){7.3}} \put(47.7,5.0){\line(-4,-1){10.7}} \put(47.7,5.0){\line(-1,3){3.2}} 
\put(55.0,12.3){\circle*{1}} \put(55.0,12.3){\line(4,1){10.0}} \put(55.0,12.3){\line(-3,1){10.5}} \put(55.0,12.3){\line(1,4){2.5}} 
\put(65.0,15.0){\circle*{1}} \put(65.0,15.0){\line(4,-1){10.0}} \put(65.0,15.0){\line(-1,1){7.5}} \put(65.0,15.0){\line(1,1){7.5}} 
\put(75.0,12.3){\circle*{1}} \put(75.0,12.3){\line(1,-1){7.3}} \put(75.0,12.3){\line(-1,4){2.5}} \put(75.0,12.3){\line(3,1){10.5}} 
\put(82.3,5.0){\circle*{1}} \put(82.3,5.0){\line(1,3){3.2}} \put(82.3,5.0){\line(4,-1){10.7}} 
\put(37.0,2.5){\circle*{1}} \put(37.0,2.5){\line(3,5){7.5}} \put(37.0,2.5){\line(-4,1){8.7}} 
\put(44.5,15.5){\circle*{1}} \put(44.5,15.5){\line(5,3){13.0}} \put(44.5,15.5){\line(-1,1){6.4}} 
\put(57.5,23.0){\circle*{1}} \put(57.5,23.0){\line(1,0){15.0}} \put(57.5,23.0){\line(-1,4){2.3}} 
\put(72.5,23.0){\circle*{1}} \put(72.5,23.0){\line(5,-3){13.0}} \put(72.5,23.0){\line(1,4){2.3}} 
\put(85.5,15.5){\circle*{1}} \put(85.5,15.5){\line(3,-5){7.5}} \put(85.5,15.5){\line(1,1){6.4}} 
\put(93.0,2.5){\circle*{1}} \put(93.0,2.5){\line(4,1){8.7}} 
\put(28.3,4.8){\circle*{1}} \put(28.3,4.8){\line(3,5){9.8}} \put(38.1,21.9){\line(-5,-1){14.7}} \put(28.3,4.8){\line(-1,3){4.9}} 
\put(38.1,21.9){\circle*{1}} \put(38.1,21.9){\line(5,3){17.0}} \put(55.2,31.7){\line(-3,1){14.2}} \put(38.1,21.9){\line(1,5){2.9}} 
\put(55.2,31.7){\circle*{1}} \put(55.2,31.7){\line(1,0){19.7}} \put(74.8,31.7){\line(-5,6){9.8}} \put(55.2,31.7){\line(5,6){9.8}} 
\put(74.8,31.7){\circle*{1}} \put(74.8,31.7){\line(5,-3){17.0}} \put(91.9,21.9){\line(-1,5){2.9}} \put(74.8,31.7){\line(3,1){14.2}} 
\put(91.9,21.9){\circle*{1}} \put(91.9,21.9){\line(3,-5){9.8}} \put(101.7,4.8){\line(1,3){4.9}} \put(91.9,21.9){\line(5,-1){14.7}} 
\put(101.7,4.8){\circle*{1}} 
\put(23.4,19.0){\circle*{1}} \put(23.4,19.0){\line(1,1){17.6}} 
\put(41.0,36.6){\circle*{1}} \put(41.0,36.6){\line(4,1){24.0}} 
\put(65.0,43.0){\circle*{1}} \put(65.0,43.0){\line(4,-1){24.0}} 
\put(89.0,36.6){\circle*{1}} \put(89.0,36.6){\line(1,-1){17.6}} 
\put(106.6,19.0){\circle*{1}}
\thinlines
\put(47.7,5.0){\line(-1,-2){2}} \put(82.3,5.0){\line(1,-2){2}}
\put(37.0,2.5){\line(1,-1){4}} \put(93.0,2.5){\line(-1,-1){4}}
\put(37.0,2.5){\line(0,-1){5}} \put(93.0,2.5){\line(0,-1){5}}
\put(28.3,4.8){\line(0,-1){5}} \put(101.7,4.8){\line(0,-1){5}}
\put(28.3,4.8){\line(-1,-1){4}} \put(101.7,4.8){\line(1,-1){4}}
\put(23.4,19.0){\line(-1,-3){2}} \put(106.6,19.0){\line(1,-3){2}}

\put(48.3,3.5){$a_2$} \put(54.5,9.7){$a_1$} \put(63.0,12.0){$a_0$} \put(71,9.7){$a_{n-1}$} \put(76,3.5){$a_{n-2}$} 
\put(35.7,5.0){$b_2$} \put(44.5,17.8){$b_1$} \put(58.5,24.5){$b_0$} \put(74.0,23.3){$b_{n-1}$} \put(86.8,14.5){$b_{n-2}$} \put(93.5,0.5){$b_{n-3}$} 
\put(25.4,5.4){$c_2$} \put(35.7,23.3){$c_1$} \put(53.6,33.6){$c_0$} \put(73.5,34){$c_{n-1}$} \put(92.3,23.3){$c_{n-2}$} \put(102.6,5.4){$c_{n-3}$} 
\put(19.7,20.0){$d_1$} \put(38.0,38.3){$d_0$} \put(63.0,45.0){$d_{n-1}$} \put(88.0,38.3){$d_{n-2}$} \put(106.3,20.0){$d_{n-3}$}

\end{picture}
    \caption{\label{fig:tn} \small The graph of convex polytope $T_n$}
\end{figure}

%% file: paper.bbl
\begin{thebibliography}{99}



\bibitemart{sla75}
{P.J. Slater}
{Leaves of trees}
{Congr. Numer}{14(37)}{1975}{549--559}

\bibitemart{har76}
{F. Harary, R. Melter}
{On the metric dimension of a graph}
{Ars Combin}{2(1)}{1976}{191--195}

\bibitemart{bacya92}
{M. Ba\'ca}
{On magic labellings of convex polytopes}
{Annals Discrete Math.}{51}{1992}{13--16}


\bibitemart{cac07}
{J. C\'aceres, C. Hernando, M. Mora, I.M. Pelayo,
M.L. Puertas, C. Seara, D.R. Wood}
{On the metric dimension of Cartesian products of graphs} 
{SIAM Journal in Discrete Mathematics}{21(2)}{2007}{423--441}

\bibitemart{gon16}
{A. Gonz{\'a}lez, C.~Hernando, M.~Mora}
{The equidistant dimension of graphs}
{Bulletin of the Malaysian Mathematical Sciences Society} {45(4)} 
 {2022} {1757--1775}

\bibitemart{np}
{A. Gispert-Fern{\'a}ndez, J.~A. Rodr{\i}guez-Vel{\'a}zquez} 
{The equidistant dimension of graphs: {NP}-completeness and the case of lexicographic product
  graphs} 
{AIMS Mathematics} {9(6)} {2024} {15325--15345}

\bibitemart{imr11b}
{M. Imran, A.Q. Baig, M.K. Shafiq, A. Semani\v{c}ov\'a-Fe\v{n}ov\v{c}ikov\'a}
{Classes of convex polytopes with constant metric dimension}
{Utilitas Math.}{90}{2013}{85--99}

\bibitemart{imr11e}
{M. Imran, U.H. Bokhary, A.Q. Baig}
{On the metric dimension of rotationally-symmetric convex polytopes}
{J. Algebra Comb. Discrete Appl.}{3(2)}{2015}{45--59}

\bibitemart{kra12a}
{J. Kratica, V. Kova\v cevi\'c-Vuj\v ci\'c, M. \v Cangalovi\'c, M. Stojanovi\'c}
{Minimal doubly resolving sets and the strong metric dimension of some convex polytopes}
{Applied Mathematics and Computation}{218}{2012}{9790--9801}

\bibitem{VV}
{J.~Kratica, M. {\v{C}}angalovi{\'c}, V. Kova{\v{c}}evi{\'c}-Vuj{\v{c}}i{\'c}}, 
{Equidistant dimension of Johnson and Kneser graphs}, arXiv:2406.17870, {2024}.

\bibitem{bal09}
{K. Balakrishnan, M. Changat, I. Peterin, S. {\v S}pacapan, P. {\v S}parl, A.R. Subhamathi.}
{Strongly distance-balanced graphs and graph products},
{European Journal of Combinatorics},
{30(5):}{1048--1053, }{2009.}

\bibitemart{imr11c}
{M. Imran, S.A.U.H. Bokhary, A.Q. Baig}
{On families of convex polytopes with constant metric dimension}
{Computers and Mathematics with Applications}{60}{2629--2638, } {2010.}

\bibitemart{cacsrom17}
{J. Casas-Roma, J. Herrera-Joancomart\'i, V. Torra} 
{A survey of graph-modiﬁcation techniques for privacy-preserving on networks} 
{Artif. Intell. Rev.} {47(3)} {341-–366, } {2017}

\bibitemart{ches11}
{S. Chester, G. Srivastava, G.} 
{Social network privacy for attribute disclosure attacks} 
In: {2011 International Conference on Advances in Social Networks Analysis and Mining}
{445–-449} {IEEE, } {2011}

\bibitemart{ches13}
{S. Chester, B.M. Kapron, G. Srivastava, S. Venkatesh} 
{Complexity of social network anonymization}
{Soc. Netw. Anal. Min.} {3(2)} {151-–166} {2013}

\bibitemart{fed08}
{T. Feder, S. Nabar, E. Terzi} 
{Anonymizing graphs} 
{CoRR arXiv:0810.5578} {2008}{}{}

\bibitemart{tru16}
{R. Trujillo-Rasua, I.G. Yero, I.G} 
{$k$-metric antidimension: a privacy measure for social graphs} 
{Inf. Sci.} {328} {403-–417} {2016}

\bibitemart{zhe11}
{E. Zheleva, L. Getoor, L.} 
{Privacy in social networks: a survey}
{Aggarwal, C.C. (ed.) Social Network Data Analytics} 
{1st edn} {277-–306} { Springer, Berlin} {2011}

\bibitemart{zho08}
{B. Zhou, J. Pei, W.S. Luk} 
{A brief survey on anonymization techniques for privacy preserving
publishing of social network data} 
{SIGKDD Explor. Newsl} {10(2)} {12-–22} {2008}

\bibitemart{mil06}
{M. Miller, M. Ba\v{c}a, J.A. MacDougall} 
{Vertex-magic total labelling of generalized Petersen graphs and convex polytopes} 
{JCMCC} {59} {89--99} {2006} 






\end{thebibliography}
